%% file: paper-elasticnet.tex
\documentclass{article}
\usepackage{spconf,amsmath,graphicx}
\usepackage[setbb,setcolon]{kmath}
\usepackage{amsmath,amssymb,amsfonts}
\usepackage{algorithmic}
\usepackage{xcolor}
\usepackage{cleveref}
\usepackage{url}

\crefname{section}{Sec.}{Sec.}
\Crefname{proposition}{Prop.}{Prop.}
\Crefname{lemma}{Lem.}{Lem.}
\Crefname{proposition}{Prop.}{Prop.}
\Crefname{corollary}{Cor.}{Cor.}
\Crefname{appendix}{App.}{App.}
\Crefname{figure}{Fig.}{Fig.}
\Crefname{algocf}{Alg.}{Alg.}
  
\usepackage[numbers,sort&compress]{natbib}

\usepackage{ifthen}
\newboolean{cameraready}
\setboolean{cameraready}{false}

\input{headings/packages}

\input{headings/shortcuts}
\input{headings/macros}

\title{Screen \& Relax: Accelerating the resolution of Elastic-net\\ by safe identification of the solution support}

\name{Théo Guyard$^{\star}$ \qquad Cédric Herzet$^{\dagger}$ \qquad Clément Elvira$^{\ddagger}$}
\address{
	$^{\star}$ Univ Rennes, INSA Rennes, CNRS, IRMAR-UMR 6625, F-35000, France \\ 
	$^{\dagger}$ INRIA Rennes-Bretagne Atlantique, Campus de Beaulieu, 35000 Rennes, France \\ 
	$^{\ddagger}$ IETR UMR CNRS 6164, CentraleSupelec Rennes Campus, 35576 Cesson Sévigné, France \\
	firstname.lastname@\{insa-rennes,inria,centralesupelec\}.fr
	\thanks{
		The research presented in this paper is reproducible.
		Code and data are available at \protect\url{https://gitlab.insa-rennes.fr/Theo.Guyard/screen-and-relax}
	}
}

\begin{document}

\maketitle
\begin{abstract}
	In this paper, we propose a procedure to accelerate the resolution of the well-known ``Elastic-Net'' problem.  
	Our procedure is based on the (partial) identification of the solution support and the reformulation of the original problem into a problem of reduced dimension. 
	The identification of the support leverages the novel concept of \emph{``safe relaxing''} where one aims to identify \emph{non-zero} coefficients of the solution. It can be viewed as a dual approach to \emph{``safe screening``} introduced in the last decade and allowing to reduce the problem dimension using the identification of \emph{zero} coefficients of the solution.
	We show numerically that combining both methodologies in a \emph{``Screen \& Relax''} strategy enables to significantly improve the tradeoff between complexity and accuracy achievable by standard resolution techniques. 
\end{abstract}

\input{sections/introduction}
\input{sections/target_problem}

\input{sections/screening}
\input{sections/relaxing}
\input{sections/ts}
\input{sections/screen_relax}
\input{sections/results}
\input{sections/conclusion}

\clearpage
\let\OLDthebibliography\thebibliography
\renewcommand\thebibliography[1]{
  \OLDthebibliography{#1}
  \setlength{\parskip}{0pt}
  \setlength{\itemsep}{6.0pt plus 0.3ex}
}

\bibliographystyle{IEEEbib}
\bibliography{references}

\end{document}

%% file: headings/packages.tex
\usepackage[utf8]{inputenc}
\usepackage[T1]{fontenc}
\usepackage{xparse}
\usepackage{pgffor}
\usepackage{ifthen}
\usepackage{xspace}
\usepackage{xcolor}
\usepackage{amsmath}
\usepackage{amssymb}
\usepackage{amsfonts}
\usepackage{amsthm}
\usepackage{mathrsfs}
\usepackage{bbm}
\usepackage{mathtools}
\usepackage{stmaryrd}
\usepackage[ruled,linesnumbered]{algorithm2e}

\SetCommentSty{mycommfont}
\usepackage{csquotes}
\usepackage{dblfloatfix}
\usepackage{comment}

%% file: headings/shortcuts.tex
\AfterEndEnvironment{definition}{\noindent\ignorespaces}

\AfterEndEnvironment{remark}{\noindent\ignorespaces}

\AfterEndEnvironment{proposition}{\noindent\ignorespaces}

\AfterEndEnvironment{lemma}{\noindent\ignorespaces}

\AfterEndEnvironment{theorem}{\noindent\ignorespaces}

\AfterEndEnvironment{corollary}{\noindent\ignorespaces}

\AfterEndEnvironment{example}{\noindent\ignorespaces}

\AfterEndEnvironment{problem}{\noindent\ignorespaces}

\foreach \x in {a,...,z}{%
	\expandafter\xdef\csname \x\endcsname{\noexpand\ensuremath{\noexpand\mathbf{\x}}}
}
\foreach \x in {a,...,z}{%
	\expandafter\xdef\csname \x rm\endcsname{\noexpand\ensuremath{\noexpand\mathrm{\x}}}
}

\foreach \x in {A,...,Z}{%
	\expandafter\xdef\csname \x\endcsname{\noexpand\ensuremath{\noexpand\mathbf{\x}}}
}
\foreach \x in {A,...,Z}{%
	\expandafter\xdef\csname \x rm\endcsname{\noexpand\ensuremath{\noexpand\mathrm{\x}}}
}
\foreach \x in {A,...,Z}{%
	\expandafter\xdef\csname \x bb\endcsname{\noexpand\ensuremath{\noexpand\mathbb{\x}}}
}
\foreach \x in {A,...,Z}{%
	\expandafter\xdef\csname \x c\endcsname{\noexpand\ensuremath{\noexpand\mathcal{\x}}}
}

\def\0{{\mathbf 0}}
\def\1{{\mathbf 1}}

\def\ie{\textit{i.e.}}
\def\eg{\textit{e.g.}}
\def\etal{\textit{et al.}}

\ifthenelse{\boolean{cameraready}}{
	\newcommand{\remCE}[1]{}

	\newcommand{\remCH}[1]{}

}{
	\newcommand{\remCE}[1]{{\noindent\color{purple}{\scriptsize[CE: #1]}}}

	\newcommand{\remCH}[1]{{\noindent\color{red}{\scriptsize[CH: #1]}}}

}

%% file: headings/macros.tex
\newcommand{\rego}{\ensuremath{\lambda}\xspace}
\newcommand{\regt}{\ensuremath{\varepsilon}\xspace}

\newcommand{\pfun}{\ensuremath{P}\xspace}
\newcommand{\pve}{\ensuremath{x}\xspace}
\newcommand{\pv}{\ensuremath{\mathbf{\pve}}\xspace}

\newcommand{\pvopt}{\ensuremath{\pv^\star}\xspace}

\newcommand{\dve}{\ensuremath{u}\xspace}
\newcommand{\dv}{\ensuremath{\mathbf{\dve}}\xspace}

\newcommand{\dvopt}{\ensuremath{\dv^\star}\xspace}

\newcommand{\obs}{\ensuremath{\y}\xspace}

\newcommand{\dic}{\ensuremath{\A}\xspace}
\newcommand{\atom}{\ensuremath{\a}\xspace}
\newcommand{\Id}{\ensuremath{\I}\xspace}

\newcommand{\pdim}{\ensuremath{n}\xspace}
\newcommand{\ddim}{\ensuremath{m}\xspace}

\newcommand{\setposopt}{\ensuremath{\mathcal{J}^\star}\xspace}

\newcommand{\setscreen}{\ensuremath{{\mathcal{I}}}\xspace}
\newcommand{\setscreenc}{\ensuremath{{\overline{\mathcal{I}}}}\xspace}
\newcommand{\setrelax}{\ensuremath{{\mathcal{J}}}\xspace}
\newcommand{\setrelaxc}{\ensuremath{{\overline{\mathcal{J}}}}\xspace}
\newcommand{\setsr}{\ensuremath{{\mathcal{K}}}\xspace}
\newcommand{\setsrc}{\ensuremath{{\overline{\mathcal{K}}}}\xspace}

\newcommand{\card}[1]{\ensuremath{\mathrm{card}(#1)}\xspace}

\newcommand{\normT}[1]{\ensuremath{\|#1\|_2}\xspace}

\newcommand{\idxscreen}{\ell}
\newcommand{\idxrelax}{\ell}

\newcommand{\safesphere}{\mathcal{S}}
\newcommand{\spherec}{\mathbf{c}}
\newcommand{\spherer}{r}

\newcommand{\matxpos}{\B}

\newcommand{\vecxpos}{\b}

\newcommand{\reducedNormTwo}{\M}
\newcommand{\obsred}{\obs_r}
\newcommand{\pvred}{\pv_r}
\newcommand{\dicred}{\dic_r}
\newcommand{\pdimred}{\pdim_r}



%% file: sections/introduction.tex

\section{Introduction} \label{sec:introduction}

Sparse decomposition aims at finding some approximation of a vector \(\obs \in \kR^{\ddim}\) as the linear combination of a few columns (dubbed \emph{atoms}) of a dictionary \(\dic=[\atom_1,\dotsc,\atom_\pdim]\in\kR^{\ddim\times\pdim}\). 
Unfortunately, identifying the sparsest decomposition of a vector according to some accuracy criterion often turns out to be a combinatorial problem~\cite[Sec.~2.3]{Foucart2013aa}. 
A standard strategy to circumvent this issue consists in approximating this ideal  decomposition as the solution of a problem of the form
\begin{equation}
	\label{eq:sparse-problem}
	\pvopt\in\kargmin_{\scriptstyle\pv\in \kR^\pdim} \ \tfrac{1}{2}\kvvbar{\obs-\dic\pv}_2^2 + \Omega(\pv)
\end{equation}
where \(\kfuncdef{\Omega}{\kR^\pdim}{\kR+}\) is some \emph{sparsity-inducing} convex regularizer. 
The common choice \(\Omega(\pv)=\rego\kvvbar{\pv}_1\) for some $\lambda>0$ leads to the well-known \textit{``Lasso''} problem and has been extensively studied in the literature \cite{tibshirani96regression,Chen_siam99}. 
Another standard choice is 
$\Omega(\pv)=\rego\kvvbar{\pv}_1 + \tfrac{\regt}{2}\kvvbar{\pv}_2^2$ for some parameters $\rego>0$, $\regt>0$. 
In this case, problem~\eqref{eq:sparse-problem} is known as \textit{``Elastic-Net''} and is popular in many applicative domains because its solution enjoys desirable statistical properties \cite{Zou:2005gi}.

Because of its clear practical interest, many contributions of the literature have proposed efficient solving procedures for \eqref{eq:sparse-problem}, see \eg, \cite{Figueiredo2007Gradient,parikh2014proximal,Jin:2009wj,boyd2011distributed}.
Of particular interest in this paper is the ``safe screening'' acceleration technique proposed by El Ghaoui \etal\ in~\cite{ghaoui2010safe}. 
Safe screening consists in performing simple tests to identify the \textit{zero} elements of the minimizers of an optimization problem. This knowledge can then be exploited to reduce the dimensionality of the problem by discarding the atoms of the dictionary weighted by zero safely, \ie, without changing the solution set. 
Over the past decade, many authors have identified safe screening as a simple procedure to significantly speed up the resolution of many optimization problems, see \eg,~\cite{Fercoq_Gramfort_Salmon15, Xiang:2017ty, icml2014c2_liuc14,wang2013lasso, Herzet16Screening,Herzet:2019fj}.

In this paper, we introduce a dual approach to safe screening, dubbed \textit{``safe relaxing''}. We focus on a specific instance of problem~\eqref{eq:sparse-problem}, namely the non-negative version of Elastic-Net. 
Our method aims at identifying the position of the \emph{non-zero} coefficients of the minimizer of this problem.  
We show that, similarly to screening, this knowledge can be exploited to safely reduce the dimensionality of the target problem and accelerate its resolution.
We use the terminology ``relaxing'' as the reduction of the problem dimensionality results from the relaxation of some constraints.

The rest of the paper is organized as follows. 
The target problem is defined in \Cref{sec:target problem}. The concepts of ``safe screening'' and ``safe relaxing'' are presented in \Cref{sec:safe screening,sec: safe relaxing}. In \Cref{sec:SR algorithm}, we combine screening and relaxing methodologies in a \textit{``Screen \& Relax''} strategy. A numerical evaluation of the proposed method is finally carried out in \Cref{sec:results}. \\

\textbf{Notations.} Boldface uppercase (\eg, $\dic$) and lowercase (\eg, $\pv$) letters respectively represent  matrices and vectors. 
$\0_\pdim$ and $\1_\pdim$ stand for the $\pdim$-dimensional all-zeros and all-ones vectors. 
$\Id$ represents the identity matrix whose dimension is usually clear from the context. 
The $i$th component of $\pv$ is denoted $\pv(i)$. 
Calligraphic letters (\eg, $\Ic$) are used to denote sets and
the notation $\overline{\Ic}$ refers to the complementary set of $\Ic$. 
We denote by $\pv_\Ic$ the restriction of $\pv$ to its elements indexed by $\Ic$ and $\dic_\Ic$ corresponds to the restriction of $\dic$ to its columns indexed by $\Ic$. Finally, for any real symmetric positive definite matrix $\reducedNormTwo$, we let $\|\pv\|^2_{\reducedNormTwo}\triangleq\ktranspose{\pv}\reducedNormTwo\pv$. Throughout this paper, we assume without loss of generality that the columns of $\dic$ are normalized to one.

%% file: sections/target_problem.tex

\section{Target Problem}\label{sec:target problem} 

We focus on the non-negative version of Elastic-Net:
\begin{align} \label{eq:primal} 
	\stepcounter{equation}
	\tag{\theequation-\(\Pc\)}
	\min_{\scriptstyle\pv \geq \0_{\pdim}} \ 
			\pfun(\pv)
			\triangleq
			\tfrac{1}{2} \kvvbar{\obs - \dic\pv}_2^2 
		    				+ \ktranspose{\boldsymbol{\rego}} \pv
		    				+ \tfrac{\regt}{2}\kvvbar{\pv}_2^2 
\end{align}
where $\boldsymbol{\rego}\in\kR+^{\pdim}$ and $\regt>0$. We note that the standard formulation of Elastic-Net can be seen as a particular case of \eqref{eq:primal} (see \eg,~\cite[Sec. 2]{xiang2016screening}). 
Since $\pfun(\cdot)$ is continuous, coercive and strongly convex,  
\eqref{eq:primal} admits a unique minimizer $\pvopt$. 
The goal of this paper is to accelerate the resolution of \eqref{eq:primal} by identifying the position of the zero and non-zero coefficients of $\pvopt$. 
Our strategy leverages the primal-dual optimality conditions described below.

The dual problem associated to \eqref{eq:primal} reads
\begin{align} 
		\label{eq:dual}
		\stepcounter{equation}
		\tag{\theequation-\(\Dc\)}
		\nonumber
		\max_{\scriptstyle\dv\in\kR^{\ddim}} 
			\tfrac{1}{2}
				\normT{\obs}^2 
				- \tfrac{1}{2}\normT{\obs - \dv}^2 
				-\tfrac{1}{2\regt}
				\|[\ktranspose{\dic}\dv - \boldsymbol{\rego}]_+\|_2^2
\end{align}
where \([\pv]_+ \triangleq \max(\0_\pdim, \pv)\) and with the maximum taken component-wise~\cite[Sec.~5.2]{dunner2016primal}.  
Similarly to \eqref{eq:primal}, the cost function in \eqref{eq:dual} is continuous, coercive and strongly concave. Problem \eqref{eq:dual} thus admits a unique maximizer $\dvopt$.~By Slater's constraint qualification, strong duality holds between \eqref{eq:primal} and \eqref{eq:dual}. 
As a consequence, a couple \((\pvopt,\dvopt)\) is a primal-dual solution of \eqref{eq:primal}-\eqref{eq:dual} if and only if
\begin{align} 
	\label{eq:optcond0}
	\dvopt &= \obs - \dic \pvopt\\
	\label{eq:optcond}
	\pvopt &= \regt^{-1} [\ktranspose{\dic}\dvopt - \boldsymbol{\rego}]_+
	.
\end{align}
See \cite[Prop.~5.1.5 and 5.3.1]{bertsekas1997nonlinear} for technical details.
In particular, letting $\setposopt \triangleq \kset{\idxrelax}{\pvopt(\idxrelax) > 0}$,  
we also easily obtain from \eqref{eq:optcond0}-\eqref{eq:optcond} that 
\begin{align}\label{eq:pvopt nonzero}
		\pvopt_{\setposopt}
		&= 
		(\ktranspose{\dic}_{\setposopt}\dic_{\setposopt} + \regt\Id)^{-1}
		(
		{\ktranspose{\dic}_{\setposopt}\obs 
		- \boldsymbol{\rego}_{\setposopt}}
		)
		.
\end{align}

%% file: sections/screening.tex

\section{Safe screening}\label{sec:safe screening}

The goal of safe screening is to identify the zero components of $\pvopt$ in order to transform 
\eqref{eq:primal} into a 
problem of reduced dimension and speed-up its resolution. 
More precisely, let 
\begin{align} \label{eq:subset zeros}
	\setscreen \subseteq  \kset{\idxscreen}{\pvopt(\idxscreen) = 0}
\end{align}
denote a subset of the zero components of $\pvopt$. 
Then, \eqref{eq:primal} is equivalent to 
	\begin{align} \label{eq:screened primal full 3} 
		\pvopt = \kargmin_{\pv\in\kR^\pdim} 
				&
				\ 
					\pfun(\pv)
				\
				\mbox{ s.t.}
				\
				\left\{
				\begin{array}{ll}
					\pv_{\setscreenc} &\geq \0_{\pdimred}\\
					\pv_{\setscreen} &= \0_{\pdim-\pdimred}
				\end{array}
				\right.
	\end{align}	
where $\pdimred \triangleq \card{\setscreenc}$.
This problem can also be rewritten more explicitly as
\begin{subequations}
	\begin{align}
			\pvopt_{\setscreenc} \label{eq: reduced problem screening 1}
			&= \kargmin_{\scriptstyle\pvred \geq \0_{\pdimred}} \ 
					\tfrac{1}{2} \kvvbar{\obs - \dicred \pvred}_2^2 
					+ \ktranspose{\boldsymbol{\rego}}_{r}
					\pvred 
					+ \tfrac{\regt}{2}\kvvbar{\pvred}_2^2 
				\\
			\pvopt_{\setscreen} \label{eq: reduced problem screening 2}
			& =	\0_{\pdim-\pdimred}	
	\end{align}	
\end{subequations}	
where $\dicred \triangleq \dic_{\setscreenc}\in\kR^{ \ddim \times \pdimred}$ and $\boldsymbol{\rego}_r\triangleq\boldsymbol{\rego}_{\setscreenc}\in\kR^{\pdimred}$. 
In the above formulation, we note that \eqref{eq: reduced problem screening 1} has the same structure as \eqref{eq:primal} but with a reduced optimization domain of dimension $\pdimred$ instead of $\pdim$. 
Hence, if $\pdimred\ll \pdim$, huge computational savings can potentially be achieved by considering the reduced formulation \eqref{eq: reduced problem screening 1} instead of \eqref{eq:primal}. 

Safe screening tests aim to identify some subset $\setscreen\subseteq\{1,\ldots,\pdim\}$ verifying \eqref{eq:subset zeros}.~The design of such tests usually leverages the optimality conditions of the problem at stake. 
As far as \eqref{eq:primal} is concerned, we have from \eqref{eq:optcond} that
\begin{align} \label{eq:ideal-screening-test bis}
	\forall \idxscreen\in \{1,\ldots,\pdim\}:\
	\ktranspose{\atom}_\idxscreen\dvopt \leq \boldsymbol{\rego}(\idxscreen) \iff \pvopt(\idxscreen)=0
	.
\end{align}
The left-hand side of the equivalence is thus a sufficient condition for $\pvopt(\idxscreen)$ to be equal to zero. 
Unfortunately, computing $\dvopt$ is usually as difficult as solving primal problem~\eqref{eq:primal} and \eqref{eq:ideal-screening-test bis} is therefore of poor practical interest.  

This difficulty can be circumvented by using \textit{``safe regions''}, that is subsets of the dual domain that are guaranteed to contain $\dvopt$.~For example, assuming that $\dvopt$ belongs to a safe spherical regions, that is 
\begin{align}\label{eq:safe sphere definition}
	\dvopt 
	\in 
	\safesphere(\spherec,\spherer)
	\triangleq 
	\kset{\dv\in\kR^\ddim}{\|\dv-\spherec\|_2\leq \spherer}
	,
\end{align}
test \eqref{eq:ideal-screening-test bis} can be relaxed as 
\begin{equation} \label{eq:screening-test B 2}
	\max_{\dv \in \safesphere(\spherec,\spherer)} \
	\ktranspose{\atom}_\idxscreen \dv
	=
	\ktranspose{\atom}_\idxscreen \spherec + \spherer
	\leq
	\boldsymbol{\rego}(\idxscreen) \implies \pvopt(\idxscreen) = 0
	.
\end{equation}
Methods to construct safe spheres have been extensively studied in the literature over the past decade, see \eg, \cite{ghaoui2010safe, Dai2012Ellipsoid, Fercoq_Gramfort_Salmon15, Xiang:2017ty, icml2014c2_liuc14,wang2015lasso, Herzet16Screening,Herzet:2019fj}.

%% file: sections/relaxing.tex

\vspace*{-1em}

\section{Safe relaxing}\label{sec: safe relaxing}

In this section, we expose our \textit{``safe relaxing''} methodology. 
In contrast to safe screening, our goal is to identify the positions of the \textit{non-zero} coefficients of $\pvopt$. We show that the identification of these components can also lead to an equivalent problem of reduced dimension. More precisely, let 
	\begin{align}\label{eq:subset nonzero}
		\setrelax \subseteq  \kset{\idxrelax}{\pvopt(\idxrelax) > 0}
	\end{align}
denote a subset of non-zero components of $\pvopt$. 
Problem~\eqref{eq:primal} can then be equivalently expressed as
	\begin{align} \label{eq:relaxed primal full 3a} 
		\pvopt = \kargmin_{\pv\in\kR^\pdim} 
				&
				\ 				
				\pfun(\pv)
				\
				\mbox{ s.t.}\
				\left\{
				\begin{array}{ll}
					\pv_{\setrelaxc} &\geq \0_{\pdimred}\\
					\pv_{\setrelax} &\in \kR^{\pdim-\pdimred}
				\end{array}
				\right.
	\end{align}
where $\pdimred \triangleq \card{\setrelaxc}$.
We note that the constraints on the elements in $\setrelax$ have been totally removed in \eqref{eq:relaxed primal full 3a}. 
This is in contrast with screening where the elements $\pvopt_\setscreen$ were set to zero. 
Similarly to screening, this relaxation allows to express \eqref{eq:primal} as a problem of reduced dimension. 

Let us first notice that the restriction of \eqref{eq:relaxed primal full 3a} to $\pvopt_\setrelaxc$ can be written as :
	\begin{align} \label{eq:relaxed primal full 3} 
		\pvopt_{\setrelaxc} = \kargmin_{\pv_{\setrelaxc}\geq \0_{\pdimred}} 
				&
				\kparen{
					\min_{\pv_\setrelax\in\kR^{\pdim-\pdimred}} 
				\pfun(\pv)\
				}
				.
	\end{align}	
Since the inner minimization in \eqref{eq:relaxed primal full 3} is a strongly-convex quadratic problem, it admits the unique optimizer
\begin{align}
	\label{eq:relaxing-relation}
	\pv_{\setrelax} &= \matxpos \pv_{\setrelaxc} + \vecxpos
\end{align}
where 
\begin{subequations}
	\begin{align} \label{eq:relaxing-relation-data}
		\matxpos &\triangleq -(\ktranspose{\dic}_{\setrelax}\dic_{\setrelax} + \regt\Id)^{-1}\ktranspose{\dic}_{\setrelax}\dic_{\setrelaxc}\\
		\vecxpos &\triangleq -(\ktranspose{\dic}_{\setrelax}\dic_{\setrelax} + \regt\Id)^{-1}\kparen{\ktranspose{\dic}_{\setrelax}\obs - \boldsymbol{\rego}_{\setrelax}}
		. \label{eq:relaxing-relation-data-B}
	\end{align}
\end{subequations}
Plugging \eqref{eq:relaxing-relation} into the cost function $\pfun(\pv)$ 
then leads to the following equivalent formulation of $\eqref{eq:primal}$:
\begin{subequations}
	\begin{align}
		\pvopt_{\setrelaxc}
		&= 
		\kargmin_{\pvred\geq \0_{\pdimred}}			
			\tfrac{1}{2} \kvvbar{\obsred - \dicred\pvred}_2^2 
			+ \ktranspose{\boldsymbol{\lambda}}_r\pvred
			+
			\tfrac{\regt}{2}
			\|\pvred\|^2_{\reducedNormTwo}
			\label{eq: relax reduced problem A}
		\\
		\pvopt_{\setrelax} 
		&= 
		\matxpos \pvopt_{\setrelaxc} + \vecxpos
		\label{eq: relax reduced problem B}
	\end{align}	
\end{subequations}	
where 
\begin{subequations}
	\begin{align}
		\dicred &\triangleq 
			\dic_{\setrelaxc} 
		    + \dic_{\setrelax} \matxpos\label{eq:update Ar}\\
		\boldsymbol{\lambda}_r &\triangleq 
		    \boldsymbol{\rego}_{\setrelaxc} 
		    + \ktranspose{\matxpos}
		    \kparen{\boldsymbol{\rego}_{\setrelax} + \regt\vecxpos}\\
		\obsred &\triangleq 
			\obs - \dic_{\setrelax}\vecxpos\\    
		\reducedNormTwo &\triangleq
			 \Id + \ktranspose{\matxpos}\matxpos \label{eq:def reducedNormTwo}
		.
	\end{align}
\end{subequations}	
Similarly to screening, the reduced problem \eqref{eq: relax reduced problem A} has the same mathematical structure as \eqref{eq:primal}. The definition of the parameters $(\dicred,\boldsymbol{\lambda}_r,\obsred,\reducedNormTwo)$ in \eqref{eq: relax reduced problem A} differs however quite significantly from those in \eqref{eq: reduced problem screening 1}.~In particular, whereas the construction of $\dicred$  only requires to remove some columns from $\dic$ in \eqref{eq: reduced problem screening 1}, it involves a matrix inversion in \eqref{eq: relax reduced problem A}. This operation introduces some complexity overhead and must therefore be performed with care as discussed in \Cref{sec:results}.

Optimality condition \eqref{eq:optcond} can be exploited to identify some subset $\setrelax$ verifying \eqref{eq:subset nonzero}. 
In particular, we have
\begin{equation} \label{eq:ideal-relaxing-test bis}
	\forall \idxscreen\in \{1,\ldots,\pdim\}:\
	\ktranspose{\atom}_\idxscreen\dvopt > \boldsymbol{\rego}(\idxscreen) \iff \pvopt(\idxscreen)>0
	.
\end{equation}
Similarly to screening, we can resort to a safe sphere \eqref{eq:safe sphere definition} to obtain a weaker, yet practical, version of \eqref{eq:ideal-relaxing-test bis}. 
This leads to the following \emph{relaxing test}:
\begin{equation} \label{eq:relaxing-test}
	\min_{\dv \in \safesphere(\spherec,\spherer)} \ \ktranspose{\atom}_\idxrelax\dv 
	=
	\ktranspose{\atom}_\idxrelax\spherec - \spherer
	>
	\boldsymbol{\rego}(\idxscreen) \implies \pvopt(\idxscreen) > 0
	.
\end{equation}

%% file: sections/ts.tex

\begin{algorithm}[t]
	\DontPrintSemicolon
	\SetKwInOut{Input}{Input}
	\Input{$\pv^{(0)}$, $\dic$, $\obs$, $\boldsymbol{\rego}$, $\regt$}
	$t \leftarrow 1$\;
	\((\setscreen,\setrelax,\setsr) \leftarrow (\emptyset,\emptyset,\emptyset)\)\;\;
	$(\dicred,\boldsymbol{\lambda}_r,\obsred,\reducedNormTwo) \leftarrow (\dic, \obs, \boldsymbol{\rego}, \Id)$\;

	\While{convergence criterion is not met}{
		\(\pv^{(t)}_{\setsrc} \leftarrow\) DescentStep({\(\pv^{(t-1)}_{\setsrc}, \dicred, \obsred, \boldsymbol{\rego}_r, \reducedNormTwo, \regt\)}) \label{line:descent step}\;
		Compute a new safe sphere $\safesphere(\spherec^{(t)},\spherer^{(t)})$\label{line:safe sphere update}\;
		Update $\setscreen$ with test~\eqref{eq:screening-test B 2}\label{line: screening}  \tcp*{Screening test}
		Update $\setrelax$ with test~\eqref{eq:relaxing-test}\label{line: relaxing} \tcp*{Relaxing test}
		\(\setsr \leftarrow \setscreen\cup\setrelax\)\;
		Update $\dicred, \obsred, \boldsymbol{\rego}_r, \reducedNormTwo$ with \eqref{eq:update Ar}-\eqref{eq:def reducedNormTwo}\label{line: problem update}\;
		$t \leftarrow t + 1$
	}
	\caption{``Screen \& Relax'' solving procedure}
	\label{algo:ts}
\end{algorithm}

%% file: sections/screen_relax.tex

\section{Screen and relax}\label{sec:SR algorithm}

The ``screening'' and ``relaxing'' procedures described in \Cref{sec:safe screening,sec: safe relaxing} can obviously be combined in a \emph{``Screen \& Relax''} strategy to benefit from the identification of \textit{both} zero and non-zero components of $\pvopt$.
More precisely, let \(\setscreen\) and \(\setrelax\) be subsets respectively verifying \eqref{eq:subset zeros} and \eqref{eq:subset nonzero} and let \(\setsr \triangleq (\setscreen \cup \setrelax)\) be the set of components of \(\pvopt\) already identified as zero or non-zero.
Applying the same reasoning as in \Cref{sec:safe screening,sec: safe relaxing}, we then obtain that \eqref{eq:primal} is equivalent to
\begin{subequations}
	\begin{align}
		\pvopt_{\setsrc} 
		&= 
		\kargmin_{\pvred\geq \0_{\pdimred}}	
			\tfrac{1}{2} \kvvbar{\obsred - \dicred\pvred}_2^2 
			+ \ktranspose{\boldsymbol{\lambda}}_r\pvred
			+
			\tfrac{\regt}{2}
			\|\pvred\|^2_{\reducedNormTwo}
			\label{eq: screen-relax reduced problem A}
		\\
		\pvopt_{\setrelax} 
		&= 
		\matxpos \pvopt_{\setsrc} + \vecxpos
		\label{eq: screen-relax reduced problem B}\\
		\pvopt_{\setscreen} 
		&= \0_{\card{\setscreen}} \label{eq: screen-relax reduced problem C}
		,
	\end{align}	
\end{subequations}		
where the parameters $(\dicred,\boldsymbol{\lambda}_r,\obsred,\reducedNormTwo)$ are defined as in \eqref{eq:relaxing-relation-data}-\eqref{eq:def reducedNormTwo} by using $\setsrc$ instead of $\setrelaxc$.~The dimension of reduced problem \eqref{eq: screen-relax reduced problem A} is equal to $\pdimred=\card{\setsrc}$ and thus benefits from the identification of both the zero and non-zero components of $\pvopt$ in its dimensionality reduction. 
 
Quite interestingly, when equality holds in \eqref{eq:subset zeros} and \eqref{eq:subset nonzero}, relations \eqref{eq: screen-relax reduced problem B}-\eqref{eq: screen-relax reduced problem C} entirely define the solution of \eqref{eq:primal}. 
In this case, \eqref{eq: screen-relax reduced problem B} reduces to \eqref{eq:pvopt nonzero}. The solution of  \eqref{eq:primal} can therefore be computed to \emph{machine-precision} via simple linear-algebra operations when all components of $\pvopt$ have either passed a screening or a relaxing test. 

%% file: sections/results.tex

\begin{figure*}[t]
	\centering
	\includegraphics[width=\linewidth]{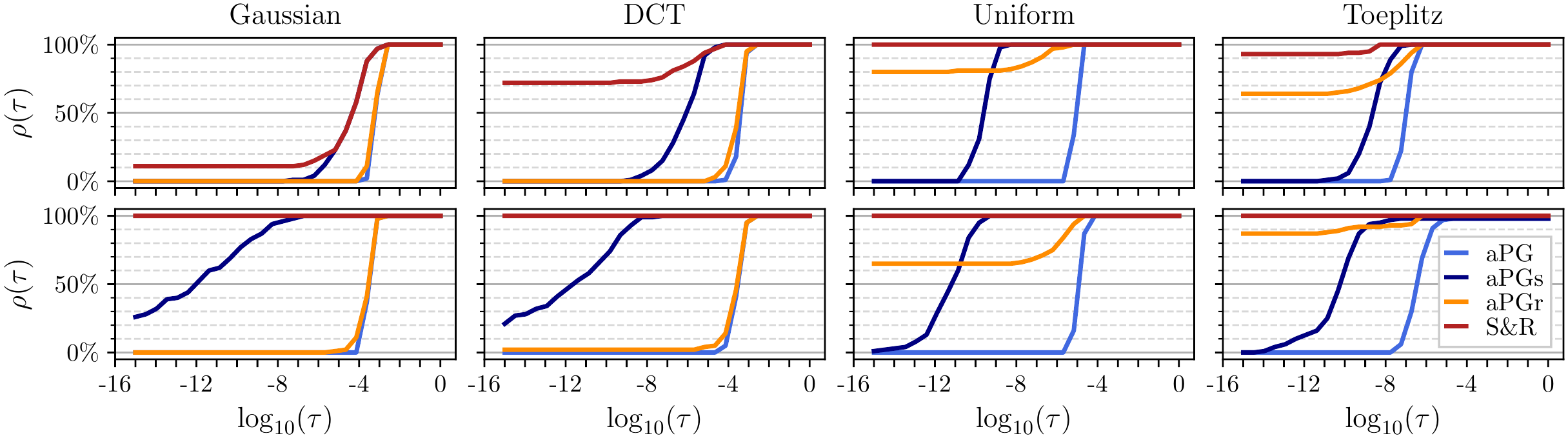}
	\caption{
		Dolan-Moré performance profiles for \((\rego,\regt)= (0.2,0.5)\rego_{\max}\) (top) and \((\rego,\regt)= (0.5,0.2)\rego_{\max}\) (bottom).
	}
	\label{fig:perf}
\end{figure*}

\section{Numerical results}\label{sec:results}

In this section, we evaluate the computational gain induced by the proposed safe relaxing strategy.
We focus on the resolution of \eqref{eq:primal} with $\boldsymbol{\rego}=\rego \1_{\pdim}$ for some $0 < \rego < \rego_{\max} \triangleq \max(\ktranspose{\dic}\obs)$. We mention that $\pvopt = \0_{\pdim}$ as soon as $\rego \geq \rego_{\max}$.

We consider the ``Screen \& Relax'' (S\&R) procedure described in \Cref{algo:ts}.
The function ``DescentStep'' in line~\ref{line:descent step} corresponds to one iteration of an accelerated proximal gradient algorithm~\cite[Sec.~4.3]{parikh2014proximal} applied to problem~\eqref{eq: screen-relax reduced problem A}. Recall that at the beginning of the solving procedure, $\setscreen = \setrelax = \emptyset$.
The evaluation of the safe sphere parameters $\spherec^{(t)}$ and $\spherer^{(t)}$ in line \ref{line:safe sphere update} follows the ``GAP'' methodology presented in \cite[Th.~6]{ndiaye2016gap}.
At each iteration, problem \eqref{eq: screen-relax reduced problem A} is updated in line \ref{line: problem update} upon the identification of additional zero or non-zero components.
We note that $\setrelax$ typically only varies by (at most) a few elements at each iteration of \Cref{algo:ts}. This behavior can be exploited to efficiently compute the inverse in \eqref{eq:relaxing-relation-data}-\eqref{eq:relaxing-relation-data-B} by using rank-one update rules~\cite{hager1989updating}.

We compare the performance of the S\&R procedure with three restricted versions of \Cref{algo:ts} : \textit{i)} no screening and no relaxing is performed (\ie, lines \ref{line: screening}-\ref{line: relaxing} are skipped); \textit{ii)} only screening is performed (\ie, line \ref{line: relaxing} is skipped); \textit{iii)} only relaxing is performed (\ie, line \ref{line: screening} is skipped). These variants will respectively be denoted ``aPG'', ``aPGs'' and ``aPGr'' in the sequel.
Both aPG and aPGs correspond to standard methodologies of the literature while aPGr and S\&R are contributions of the present paper.

We use ``Dolan-Moré'' performance profiles~\cite{dolan2002benchmarking} to assess the performance of these four methods. Our results are gathered in \Cref{fig:perf}.~To generate each curve, we run a solving method with a \textit{given} computational budget on 100 different instances of problem~\eqref{eq:primal}.~The curve corresponds to the percentage $\rho(\tau)$ of problem instances for which the solving strategy achieves a duality gap~\cite{dunner2016primal} lower than $\tau$.

To generate problem data, we consider the four following setups:
the elements of $\dic$ are i.i.d. realizations of
\textit{i)} a standard normal distribution
or \textit{ii)} a uniform law on \(\kintervcc{0}{1}\); \textit{iii)} the rows of $\dic$ are randomly-sampled from a DCT matrix~\cite{ahmed1974discrete}; \textit{iv)} $\dic$ has a Toeplitz structure~\cite{gray2006toeplitz} with shifted versions of a Gaussian curve.
In all setups, the columns of \(\dic\) are normalized to one.
The observation \(\obs\) is drawn according to a uniform
distribution on the \(\ddim\)-dimensional sphere for ``Gaussian'' and ``DCT'' dictionaries and is restricted to the positive orthant for ``Uniform'' and ``Toeplitz'' dictionaries.
We set $\ddim=100$, $\pdim=300$ and consider the following choices for the regularization parameters : \((\rego,\regt)= (0.2,0.5)\times\rego_{\max}\) or \((\rego,\regt)= (0.5,0.2)\times\rego_{\max}\). Each problem instance is solved with a budget of \(2\times10^6\) FLOPs (the number of floating-point operations) for ``Gaussian'' and ``DCT'' dictionaries, and \(2\times10^7\) FLOPs for ``Uniform'' and ``Toeplitz'' dictionaries. 
The difference in the FLOPs budgets stems from the bad conditioning of the ``Uniform'' and ``Toeplitz'' dictionaries which leads to slower convergence of standard numerical solvers.

As far as our simulation setups are concerned, we notice that safe relaxing enables us to significantly improve the performance.~Safe relaxing alone (aPGr) proves to be of particular interest for dictionaries with highly-correlated atoms (\eg, ``Uniform'' or ``Toeplitz'').
A careful study of our simulation results led us to the conclusion that this behavior is due to an improvement of the problem conditioning when moving from problem \eqref{eq:primal} to \eqref{eq: screen-relax reduced problem A} and therefore of the convergence rate of the proximal gradient algorithm.
The combination of screening and relaxing significantly outperforms all the other methods. We notice that S\&R attains machine precision ($\tau=10^{-16}$) for a large proportion of problem instances in most setups. This can be explained by the behavior emphasized in \Cref{sec:SR algorithm}: when all the zero and non-zero elements of $\pvopt$ are identified, the minimizer can be explicitly computed from \eqref{eq: screen-relax reduced problem B}-\eqref{eq: screen-relax reduced problem C} with simple linear operations. Now, perfect identification of zero and non-zero elements of $\pvopt$ always occurs after a finite number of iterations when the GAP methodology is used to construct the safe sphere in tests
\eqref{eq:screening-test B 2} and \eqref{eq:relaxing-test}
since $\spherec^{(t)}\rightarrow\dvopt$ and $\spherer^{(t)}\rightarrow 0$ as $t\rightarrow\infty$.

%% file: sections/conclusion.tex

\section{Conclusion}

In this paper, we proposed a new safe relaxing methodology to detect the position of non-zero components in the solution of the Elastic-Net problem.
We showed how to leverage this knowledge to reduce the dimension of the optimization problem, enabling potential computational gains in the resolution.
Numerical simulations show the interest of the method, especially when combined with safe screening.